\documentclass[10pt]{article}
\usepackage{amssymb}
\usepackage{latexsym}
\usepackage{exscale}

\usepackage{pgf,tikz}

\usetikzlibrary{trees}

\usepackage{amssymb, amsmath}
\usepackage{enumerate}
\usepackage{enumerate}
\usepackage[top=1in,left=1in,right=1in,bottom=1in]{geometry}

\usepackage{xcolor}
\usepackage{soul}

\usepackage{cancel}

\def \beq{\begin{equation}}
\def \eeq{\end{equation}}
\def \lab{\label}
\renewcommand{\rq}[1]{(\ref{#1})}

\newtheorem{thm}{Theorem}

\newcommand{\pa }{\partial }
\newcommand{\patt }{\partial_{tt} }
\newcommand{\paxx }{\partial_{xx}}

\newcommand{\al}{\alpha }
\newcommand{\be}{\beta }

\newcommand{\f}{\varphi }

\newcommand{\mH}{{\mathcal H}}

\newcommand{\om}{{\omega}}
\newcommand{\Om}{{\Omega}}
\newcommand{\omn}{{\omega_n}}
\newcommand{\ka}{{\kappa}}

\def\<{\langle} \def\>{\rangle}

\begin{document}

\title{ Exact Controllability for the Wave Equation on a Graph with Cycle and Delta-Prime Vertex Conditions}

 \maketitle

\noindent{\bf Sergei Avdonin}, Department of Mathematics and Statistics, University of Alaska Fairbanks, Fairbanks,
AK 99775, U.S.A.; s.avdonin@alaska.edu;  \\ and Moscow Center for Fundamental and Applied Mathematics,  Moscow 119991, Russia

\

\noindent{\bf Julian Edward}, Department of Mathematics and Statistics, Florida International University, Miami,
FL 33199, U.S.A.; edwardj@fiu.edu; 

\

\noindent{\bf G\"unter Leugering}, Departments of Mathematics and Data Science Friedrich-Alexander-University Erlangen-N\"urnberg,
Martenstr. 3, 91058 Erlangen, Germany; 	guenter.leugering@fau.de.
\

\ 

{\bf Abstract.}
Exact controllability for the wave equation on a metric graph consisting of a cycle and two attached edges is proven. One boundary and one internal control are used. At the internal vertices, delta-prime conditions are satisfied. As a second example, we examine a tripod controlled at the root and the junction, while the leaves are fixed. These examples are key to understanding controllability properties in general metric graphs.

\vskip1cm 

\section{Introduction}
Control problems  for linear wave and beam equations on metric, in particular, tree-like graphs, have been heavily researched  because of their many applications in science and engineering, see \cite{AE}, \cite{DZ}, \cite{LLS1},\cite{LLS2}, \cite{LLS3}, \cite{Z} and references therein. The theory has been extended to quasilinear wave  and beam equations in recent years, see e.g. \cite{Li} and the references therein for the wave equation, as well as \cite{Gu} for nonlinear beam equations on networks. However, until very recently, very few positive results were obtained for graphs with cycles.   {In fact, suppose controls are exerted at simple nodes only. Then linear wave equations on graphs containing cycles, or even on trees where two or more Dirichlet nodes are uncontrolled, are not  exactly controllable, and  not always even approximately controllable.} Also, applying controls at multiple nodes  is not possible as long as the classical transmission conditions, i.e. the continuity condition and Kichhoff's condition across the joints,   stay intact. This observation suggests we relax these conditions in one way or another in order to gain a stronger access to the system dynamics by the way of controls. We can distinguish two types of transmission conditions at multiple nodes which have come to be known as \emph{delta} and \emph{delta-prime} conditions in the context of electrical engineering, nanotechnology and quantum graphs. In principle, there are two ways to approach this concept: one is to relax the continuity condition of states (delta) or derivatives  (delta-prime) by inserting a control between two or more incident edges, and the second one aims at relaxing the balance of forces or states, respectively, (typically referred to as  Kirchhoff conditions). 
Significant progress in this direction recently appeared in \cite{AZ2}, where an exact controllability result for general graphs was formulated and proven. In that work, the authors assume
a mix of interior and boundary controls, and at the uncontrolled interior vertices, standard (Kirchhoff--Neumann) conditions are assumed to hold. For a related result, see \cite{AEZ}. We also mention another concept of controllability, the so-called nodal profile controllability, where at a given node the state and the velocity profile is exactly matched by boundary controls. This concept is weaker than that of exact controllability as discussed in this article, but admits cycles. See \cite{LRW} for further information and references.

In the present paper, we focus on two key elementary graphs in order to understand the controllability properties for problems on general metric graphs. These are a simple cycle with attached edges and a tripod with two Dirichlet nodes, the former being more complicated. 

\section{Problem statement}
To be more specific, we consider the graph shown in Figure 1 with  delta-prime vertex conditions holding at
the uncontrolled interior vertices. Delta-prime conditions arise naturally in problems in electrical engineering \cite{AAA} and nanotechnology \cite{Ex}. We are unaware of any control theory papers on graphs with delta-prime conditions at internal vertices, although the paper \cite{AE5} discussed an inverse problem for a quantum tree 
with delta prime conditions. If we release the junction conditions at node $v_3$, keeping only edges 3 and 4 connected and pose instead a Dirichlet condition there for edge 2, then Figure 2 emerges.
Without changing the analysis significantly, we may also replace edge 2 by two serially connected edges such that the circle becomes a triangle. In this way, we can view our two examples as building blocks for the understanding of  the controllability properties for wave equations on general metric graphs. A theory on exact controllability covering the general case, however, is still lacking and subject to further investigations.

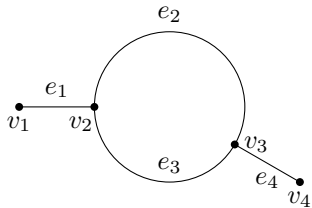
\begin{figure}
 \begin{center}
 	\begin{tikzpicture}
 	\begin{scope}[xshift=0cm]
 	\draw (0,0) circle (1cm);
 	\draw[] (-1,0) node[below] {$v_2\quad$} --node[above] {$e_1$}  (-2,0) node[below] {$v_1$}; 
 	\draw[] (330:1cm) node[right] {$v_3$} -- node[below ] {$e_4$} (330:2cm)  node[below] {$v_4$};
 	\draw[] (0,-1) node[above,black] {$e_3$};
 	\draw[] (0,1) node[above] {$e_2$};
 	\fill (-1,0) circle (1.5pt)
 	(-2,0) circle (1.5pt)
 	(-30:1cm) circle (1.5pt)
 	(-30:2cm) circle (1.5pt);;          
 	\end{scope}
 	\end{tikzpicture}
 \end{center}

		\caption{A ring with two attached edges}

\end{figure}

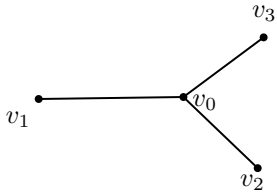
\begin{figure}	
	\begin{center}
		\tikzset{every picture/.style={line width=0.75pt}} 
		
		\begin{tikzpicture}[x=0.75pt,y=0.75pt,yscale=-1,xscale=1]

		\draw    (283.85,78.67) -- (356.85,77.67) ;
		\draw [shift={(356.85,77.67)}, rotate = 358.92] [color={rgb, 255:red, 0; green, 0; blue, 0 }  ][fill={rgb, 255:red, 0; green, 0; blue, 0 }  ][line width=0.75]      (0, 0) circle [x radius= 1.5, y radius= 1.5]   ;
		\draw    (356.85,77.67) -- (394.31,113.57) ;
		\draw [shift={(394.31,113.57)}, rotate = 43.33] [color={rgb, 255:red, 0; green, 0; blue, 0 }  ][fill={rgb, 255:red, 0; green, 0; blue, 0 }  ][line width=0.75]      (0, 0) circle [x radius= 1.5, y radius= 1.5]   ;
		\draw    (356.85,77.67) -- (397.31,47.57) ;
		\draw [shift={(397.31,47.57)}, rotate = 326.58] [color={rgb, 255:red, 0; green, 0; blue, 0 }  ][fill={rgb, 255:red, 0; green, 0; blue, 0 }  ][line width=0.75]      (0, 0) circle [x radius= 1.5, y radius= 1.5]   ;
		\draw [shift={(283.85,78.67)}, rotate = 326.58] [color={rgb, 255:red, 0; green, 0; blue, 0 }  ][fill={rgb, 255:red, 0; green, 0; blue, 0 }  ][line width=0.75]      (0, 0) circle [x radius= 1.5, y radius= 1.5]   ;

		\draw (360,75) node [anchor=north west][inner sep=0.75pt]   [align=left]{$v_0$} ;
		\draw (384.04,115.50) node [anchor=north west][inner sep=0.75pt]   [align=left] {$v_2$};
		\draw (265.85,83.50) node [anchor=north west][inner sep=0.75pt]   [align=left] {$v_1$};	
		\draw (390,30) node [anchor=north west][inner sep=0.75pt]   [align=left] {$v_3$};

		\end{tikzpicture}
		\caption{Star with coordinate system: $e_j$ identified with $(0,\ell_j)$}
	\end{center}	
\end{figure}

We begin with the graph in Figure 1, so
our graph $\Om=\{V,E\}$ consists four vertices, $V=\{v_j,\,j=1,\ldots,4\}$, and four edges,  $E=\{e_j,\,j=1,\ldots,4\}$. Each edge is identified with an interval of the real axis with positive length, $e_1(v_1,v_2)$ with $(0,l_1),$ $e_2(v_2,v_3)$ with $(0,l_2),$ $e_3(v_2,v_3)$ with $(0,l_3),$ and $e_4(v_3,v_4)$ with $(0,l_4).$
In what follows, we denote by $\phi_j$ the restriction of function $\phi$ to the edge $e_j$.
We consider the initial boundary value problem on $[\Om \setminus V ] \times [0,T]$ with one boundary control and one internal control:
\beq \lab{eqq1}
\patt u - \paxx u + q(x) u =0,
\eeq
\beq \label{ic1}
u(\cdot,0)=u_t(\cdot,0)=0,
\eeq 
\beq \lab{bcc1} 
  u_1(v_1,t)=f_1(t),
 \eeq
 \beq \lab{inc1}
\pa u_2(v_2,t)- \pa u_1(v_2,t)=f_2(t), \ \pa u_3(v_2,t)-\pa u_1(v_2,t)=0, \ \ \sum_{j=1}^3  u_j(v_2,t)=0.   
\eeq
\beq \lab{inc2}
\pa u_2(v_3,t)=\pa u_3(v_3,t)=\pa u_4(v_3,t), \ \ \sum_{j=2}^4  u_j(v_3,t)=0. 
\eeq
\beq \lab{inc4}
 u_4(v_4,t)=0. 
\eeq
Here $f_1 \in H^1(0,T),  \; f_1(0)=0, \; f_2 \in L^2(0,T),$  $q_j$ is a continuous function for all $j$, and
$\partial u_j$  denotes the derivative of $u$ at a vertex, along the incident edge $e_j$ in the direction outward from the vertex.

 We denote the space $L^2(\Omega) $ by $ \mH $ and the space of functions  on $\Omega$ whose restriction to 
each edge $e$ is in 
 $H^1(e)$ by $\mH^1, $ with norm 
 $\| \phi\|^2_{\mH^1}=\sum_{j=1}^3\| \phi '\|_{L^2(e_j)}^2+\| \phi \|_{L^2(e_j)}^2.$ 
Functions $y \in \mH^1$ such that 
 $$ y_1(v_1)=y_4(v_4)=0, \ \ \sum_{j=1}^3  y_j(v_2)=0,  \ \ \sum_{j=2}^4  y_j(v_3)=0$$
 we denote  
  by ${\mH}^1_0.$  
 For presentational simplicity, we will assume
 $$l_2\geq l_3,$$
 leaving the simple  adaptation  to the other case to the reader.
\begin{thm}\label{thm1}
Let $T_*=\max \, \{ \, l_1+l_2, \, l_3+l_4 \, \}$  and let $T\geq 2T_*$.
Let $(\phi^1,\phi^2) \in \mH_0^1 \times \mH$.
 There exist controls ${\bf f}:= (f_1,f_2)\in  H^{1}_0(0,T) \times  L^2(0,T)$ such that  the solution $u^{\bf f}$ to System \rq{eqq1}-\rq{inc4} satisfies
\beq
u^{{\bf f}}(x,T)=\phi^1(x),\ u_t^{{\bf f}}(\cdot,T)=\phi^2(\cdot),\label{ec1}
\eeq
and there exists a constant $C$ that depends only on $q,l_j,j=1,\ldots,4,$
such that 
\beq 
\| f_1\|_{H^1(0,2T)}+\| f_2\|_{L^2(0,2T)}\leq C(
\| \phi^1\|_{\mH^1}+\| \phi^2\|_{\mH}).\label{ce11}
\eeq
\end{thm}
The existence of ${\bf f}$ satisfying \rq{ec1}, \rq{ce11} is called exact controllability.

Key ingredients in the proof of Theorem 1 are shape and velocity control results that can be formulated as follows.
\begin{thm}\label{thm3}

\

a) (Shape control) Let $T\geq T_*$. Let $\phi^1 \in\mH^1_0$.
 There exist controls ${\bf f}:= (f_1,f_2)\in H_0^1(0,T)\times L^2(0,T)$ such that 
$$u^{{\bf f}}(x,T)=\phi_1(x),$$
and there exists a constant $C$ that depends only on  $q,l_j,j=1,\ldots,4,$
such that 
\beq 
\| f_1\|_{L^2(0,T)}+\| f_2\|_{H^1(0,T)}\leq C(
\| \phi_1\|_{\mH^1}).\label{ce3}
\eeq
\

b) (Velocity control). Let $T\geq T_*$. Let $\phi^2 \in \mH$.
 There exist controls ${\bf f}:= (f_1,f_2)\in H^1(0,T) \times L^2(0,T) $ such that 
$$u_t^{{\bf f}}(x,T)=\phi^2(x),$$
and there exists a constant $C$ that depends only on  $q,l_j,j=1,\ldots,4,$
such that
\beq 
\| f_1\|_{L^2(0,T)}+\| f_2\|_{H^1(0,T)}\leq C(
\| \phi^2\|_{\mH}).\label{ce4}
\eeq

\end{thm}
\vskip3mm

Theorem \ref{thm1} follows from Theorem \ref{thm3} by a moment argument.
This argument was successfully used in proving exact controllability of hybrid systems and
of the wave equation on graphs with the Kirchhoff--Neumann matching conditions \cite{AE}, \cite{AEZ}, \cite{AZ}, \cite{AZ2}. 

We remark that compared with \cite{AZ2}, this paper has a slightly sharper control time.

This paper is organized as follows. In Section 3.1, we give a solution representation of the wave equation on the interval for various boundary conditions. These representations are used in the proof of shape and velocity controllability. In Section 3.2, we prove shape controllability for a related initial boundary value problem on a 3-point star as in Figure 2. The main purpose of this section is to present the ideas of our shape control argument in a technically simpler setting, but the result might also be of independent interest. In Section 3.3, we then prove the shape controllability result announced in Theorem \ref{thm3}. In Section 3.4, we  prove Theorem \ref{thm1}. 
\section{Proofs.}
\subsection{Forward problem for the interval.}

We begin by recalling solution for the forward problem on the interval.
\begin{eqnarray}
\patt {v}-\paxx {v}+q(x){v}& = & 0,\ 0<x<l,\ t\in (0,T),\label{IBVP9a} \\
{v}(x,0)={v}_t(x,0) &= & 0,  \ 0<x<l, \\
{v}_x(0,t) & = & g(t),  \\
{v}_x(l ,t) &=  & -\pa v(l,t)=-h(t), \ t>0.\label{IBVP9b}
\end{eqnarray}
For any function $p(t)$, we denote by $P(t)$ by the formula: $P(t)=-\int_0^tp(s)ds$.
In this case, see \cite{AZ}, the solution to the system \rq{IBVP9a}-\rq{IBVP9b} is given by 
\begin{eqnarray}
{v}_{NN}^{(g,h)}(x,t)& =& 
 \sum_{n\geq 0}\left ( G(t-2nl-x)+\int_{2nl+x}^tw_+(2nl+x,s)G(t-s)ds\right )\nonumber \\
& &
+\sum_{n\geq 1}\left (G(t-2nl+x)+\int_{2nl-x}^tw_-(2nl-x,s)G(t-s)ds \right )
\nonumber \\
&&
+\sum_{n\geq 1} \left ( H(t-(2n-1)l+x) +\int_{(2n-l)\ell-x}^t w_-( (2n-1)l-x,s) H(t-s) \, ds\right ) \nonumber \\
&&
+\sum_{n\geq 1} \left ( H(t-(2n-1)l-x) +\int_{(2n-l)\ell+x}^t w_+( (2n-1)l+x,s) H(t-s) \, ds\right ), \label{nn}
\end{eqnarray}
where $w_{\pm}$ are continuous functions with $w_x\in L^2,$ and we set $G(t)=H(t)=0$
for $t<0.$

Similarly, the solution to 
\begin{eqnarray}
\patt{v}-\paxx {v}+q(x){v}& = & 0,\ 0<x<l,\ t\in (0,T),\label{IBVP9a} \\
{v}(x,0)={v}_t(x,0) &= & 0,  \ 0<x<l, \\
{v}_x(0,t) & = & g(t),  \\
{v}(l ,t) &=  & f(t), \ t>0,\label{IBVP9b}
\end{eqnarray}
is given by 
\begin{eqnarray}
{v}_{ND}^{(g,f)}(x,t)& =& 
 \sum_{n\geq 0}\mbox{ {$(-1)^n$}}\left ( G(t-2nl-x)+\int_{2nl+x}^tw_+(2nl+x,s)G(t-s)ds\right )\nonumber \\
& &+
\sum_{n\geq 1} \mbox{ {$(-1)^n$}}\left(G(t-2nl+x)+\int_{2nl-x}^tw_-(2nl-x,s)G(t-s)ds \right ),\nonumber\\
& & 
+\sum_{n\geq 1}(-1)^{n-1}\left ( f(t-(2n-1)l+x)+\int_{(2n-1)l-x}^tk_-((2n-1)l-x,s)
 f_1(t-s)ds\right )\nonumber \\
& &+
\sum_{n\geq 1}(-1)^{n-1}\left (f(t-(2n-1)l-x)+\int_{(2n-1)l+x}^tk_+((2n-1)l+x,s)f_1(t-s)ds \right ),\nonumber\\
&&
\label{nd}
\end{eqnarray}
where the properties of $k_{\pm}$ are analogous to those of $w_{\pm}$, and we set $f(t)=0$ for $t<0.$

\subsection{Shape control for three-point star}

In this section, we demonstrate our shape control argument on a simpler example: a three edge star.
In what follows, we denote by $v_0$ the central vertex, and by $v_j$ the boundary vertex for $e_j$, with $j=1,2,3, $ see Figure 2.

We consider the initial boundary value problem  with one boundary control and one internal control:
\beq \lab{eqq1b'}
\patt u - \paxx u + q(x) u =0,\ x\in e_j, t>0
\eeq
\beq \label{ic1b'}
u(\cdot,0)=u_t(\cdot,0)=0,
\eeq 
\beq \lab{bcc1b'} 
  u_1(v_1,t)=f_1(t),\  u_2(v_2,t)= u_3(v_3,t)=0.
 \eeq
 \beq \lab{inc1b'}
\pa u_2(v_0,t)- \pa u_1(v_0,t)=f_2(t), \ \pa u_3(v_0,t)-\pa u_1(v_0,t)=0, \ \ \sum_{j=1}^3  u_j(v_0,t)=0.   
\eeq
For simplicity of exposition, we assume $l_2\geq l_3$ in this section.
Let $\mH^1_0$ be a space of functions $y$ on this graph whose restriction $y_j$ to each $e_j$ belongs to $H^1(e_j)$ and $$y_1(v_1)=y_1(v_1)=y_1(v_1)=0, \ \ \sum_{j=1}^3 y_j(v_0)=0.$$
\begin{thm}
 (Shape control) Let $T_*=\max (l_2,l_1+l_3).$ Let $T\geq T_*$. Let $\phi \in \mH^1_0$.
 There exist controls ${\bf f}:= (f_1,f_2)\in H_0^1(0,T)\times L^2(0,T)$ such that 
$$u^{{\bf f}}(x,T)=\phi (x),$$
and there exists a constant $C$ that depends only on  $q,l_j,j=1,\ldots,3,$
such that 
\beq 
\| f_1\|_{L^2(0,T)}+\| f_2\|_{H^1(0,T)}\leq C(
\| \phi\|_{\mH^1}).\label{ce3}
\eeq
\

\end{thm}
Proof: 

Set $\pa g_j(t)=u_j(v_0,t).$
Applying \rq{inc1b'}, we get the equations:
\begin{eqnarray}
g_1(t) +f_2(t)& = & g_2(t),\label{g1g2b}\\
g_1(t)& = & g_3(t).\label{g1g3b}
\end{eqnarray}
For readability of what follows, we now set $q=0$, so that the integral terms in \rq{nd} vanish. The necessary modifications in the general case will be given at the end of section \ref{sc}.
Using $\sum_1^3\partial_tu_j(0,t)=0$, we get the equation:
\begin{eqnarray}
0& = &
g_1(t)+2\sum_{n\geq 1}\mbox{ {$(-1)^n$}}g_1(t-2nl_1) +2\sum_{n\geq 1}(-1)^{n}f_1'(t-(2n-1)l_1)\nonumber\\
& & +g_2(t)+2\sum_{n\geq 1}\mbox{ {$(-1)^n$}}g_2(t-2nl_2)\nonumber\\
&&+g_3(t)+2\sum_{n\geq 1}\mbox{ {$(-1)^n$}}g_3(t-2nl_3).\label{key1xxx}
\end{eqnarray}
Combining these with \rq{g1g2b}, \rq{g1g3b}, we get
$$
0=3g_1(t)+f_2(t)+2\sum_{n\geq 1}\mbox{ {$(-1)^n$}}g_1(t-2nl_1)
+2\sum_{n\geq 1}\mbox{ {$(-1)^n$}}(g_1(t-2nl_2)+f_2(t-2nl_2))+2\sum_{n\geq 1}\mbox{ {$(-1)^n$}}g_1(t-2nl_3)
$$
\beq
+2\sum_{n\geq 1}(-1)^{n}f_1'(t-(2n-1)l_1).\label{key2b}
\eeq 

The construction of the controls is a two step process using the equations above.

{\bf Step 1.}

In this step, we determine 
 $f_1^1$ supported in $(T_*-l_1-l_2,T_*-l_1)$, together with $f_2(t)$, supported on $(T_*-l_2,T_*),$
so that 
$$u^{(f_1^1,f_2)}_j(x,T_*)=\phi_j,\ j=2,3.$$

Since $q=0$, the integral kernels in  \rq{nd} vanish, so  we can easily find 
$g_j\in L^2(T_*-l_j,T_*)$ such that for $j=2,3$,
\beq
u_j(x,T_*)=v_{ND}^{(g_j,0)}(x,T_*)=\phi_j(x).\label{int}
\eeq 
We then set $g_3(t)=0$ for $t\leq  T_*-l_3,$ and $g_2(t)=0$ for $t\leq  T_*-l_2.$
Then $f_2$ is found immediately for all $t$  by applying  \rq{g1g2b},\rq{g1g3b}, and clearly  $f_2$ has its support in $[T_*-l_2,T_*]$.

We now solve for $(f_1^1)'$. Denote by $\al (t)$
various functions that are now known for  $t\in [T_*-l_2,T_*]$, which include $g_1$ and $f_2$. Thus
  \rq{key1xxx} simplifies to
$$
\sum_{n\geq 1} (-1)^n(f_1^1)'(t-(2n-1)l_1)=\al (t), \ t\in (T_*-l_2,T_*),
$$
hence 
\beq
\sum_{n\geq 0} (-1)^n(f_1^1)'(t-2nl_1)=\al (t), \ t\in (0,T_*-l_1).\label{key4'b}
\eeq 
We solve for $f_1'$ iteratively. For $t<2l_1$ and $n>0$, $(f_1^1)'(t-2nl_1)=0$ so \rq{key4'b} simplifies to 
\beq \label{known}
(f_1^1)'=\al (t);
\eeq
thus we have solved for $(f_1^1)'(t)$ for $t<2l_1.$ Now suppose we have solved for $(f_1^1)'(t)$ for $t<2n_0l_1.$ Then for $t < 2(n_0+1)l_1$ and $n>0$, $(f_1^1)'(t-2nl_1)$ is known, so again we have \rq{known} holding, and we thus solve 
for $(f_1^1)'(t)$ for $t < 2(n_0+1)l_1$. This process can be iterated until $t> T_*-l_1$, where we set $(f_1^1)'(t)=0$. Setting $f_1^1(0)=0$, this
uniquely determines $f_1^1$ on $(0,T_*-l_1).$
We now use any extension of $f_1^1$ to $H_0^1(0,T_*).$

{\bf Step 2.}
We have defined $f_2(t)$ for all $t$.
It remains to find $f_1(t)$.
Let $\tilde{\phi}=\phi -u^{(f_1^1,f_2)}(*,T_*), $ so $\tilde{\phi}_j=0$ for $j=2,3.$
By continuity, we have $\tilde{\phi}(v_0)=0$, and also $\tilde{\phi}(v_1)=0$.
Then by  \rq{nd}, we find $f_1^2\in H^1_0(0,T_*)$, with support in $[T_*-l_1,T_*]$,  so that $v_{ND}^{(0,f_1^2)}(x,T_*)=\tilde{\phi}_1(x),$ and hence
$u=u^{(f_1^2,0)}$ will solve $u_1(x,T_*)=\tilde{\phi}_1(x).$ We then  set $f_1=f_1^1+f_1^2$, and the shape control problem is solved. $\Box$

One can similarly prove velocity controllability, in time $T_*$, and then by using the arguments of the Section \ref{proof}, we get
\begin{thm}\label{thm4}
Let $T_*=\max \, \{ \, l_2, \, l_1+l_3 \, \}$  and let $T\geq 2T_*$.
Let $(\phi^1,\phi^2) \in \mH_0^1 \times \mH$.
 There exist controls ${\bf f}:= (f_1,f_2)\in  H^{1}_0(0,T) \times  L^2(0,T)$ such that  the solution $u^{\bf f}$ to System \rq{eqq1}-\rq{inc4} satisfies
\beq
u^{{\bf f}}(x,T)=\phi^1(x),\ u_t^{{\bf f}}(\cdot,T)=\phi^2(\cdot),\label{ec}
\eeq
and there exists a constant $C$ that depends only on $q,l_j,j=1,\ldots,4,$
such that 
\beq 
\| f_1\|_{H^1(0,2T)}+\| f_2\|_{L^2(0,2T)}\leq C(
\| \phi^1\|_{\mH^1}+\| \phi^2\|_{\mH}).\label{ce1}
\eeq
\end{thm}

\subsection{Shape control for graph with cycle}\label{sc}
In this section we prove part a) of Theorem \ref{thm3}. The proof for part b) is velocity controllability is similar, and is left to the reader.
We begin by using \rq{inc1},\rq{inc2}, to derive some equations, \rq{g1g2}-\rq{key3}, that will be used in constructing the shape control.

Now let $u$ be the solution to the system \rq{eqq1}-\rq{inc4}. We adopt the following notation:
\beq
g_j(t)=\partial u_j(v_2,t), \ j=1,2,3,\mbox{ and } h_k(t)=\partial u_k(v_3,t),\ k=2,3,4.
\label{gh}.
\eeq 
Then for each $j$, we can use \rq{nn} or \rq{nd} to represent $u_j$. In particular,
\beq
u_1(x,t)=v^{(g_1,f_1)}_{ND}(x,t);\ u_j(x,t)=v_{NN}^{(g_j,h_j)}(x,t)\mbox{ for }j=2,3;\ 
u_4(x,t)=v_{ND}^{h_4,0}(x,t).\label{not}
\eeq
Applying \rq{inc1}, we get the equations:
\begin{eqnarray}
g_1(t) +f_2(t)& = & g_2(t),\label{g1g2}\\
g_1(t)& = & g_3(t).\label{g1g3}
\end{eqnarray}
For readability of what follows, we now set $q=0$, so that the integral terms in \rq{nn} and \rq{nd} vanish. The necessary modifications in the general case will be given at the end of section \ref{sc}.
Using $\sum_1^3\partial_tu_j(0,t)=0$, we get the equation:
\begin{eqnarray}
0& = &
g_1(t)+2\sum_{n\geq 1}\mbox{ {$(-1)^n$}}g_1(t-2nl_1) +2\sum_{n\geq 1}(-1)^{n-1}f_1'(t-(2n-1)l_1)\nonumber\\
& & +g_2(t)+2\sum_{n\geq 1}g_2(t-2nl_2)+
2\sum_{n\geq 1}h_2(t-(2n-1)l_2)\nonumber\\
&&+g_3(t)+2\sum_{n\geq 1}g_3(t-2nl_3)
+ 2\sum_{n\geq 1}h_3(t-(2n-1)l_3).\label{key1}
\end{eqnarray}
Combining these and using $h_2(t)=h_3(t)=h_4(t)$  , we get
$$
0=3g_1(t)+f_2(t)+2\sum_{n\geq 1}g_1(t-2nl_1)
+2\sum_{n\geq 1}(g_1(t-2nl_2)+f_2(t-2nl_2))+2\sum_{n\geq 1}g_1(t-2nl_3)
$$
\beq
+2\sum_{n\geq 1}(-1)^{n-1}f_1'(t-(2n-1)l_1)+2\sum_{n\geq 1}h_4(t-(2n-1)l_2)+ 2\sum_{n\geq 1}h_4(t-(2n-1)l_3),\forall t.\label{key2}
\eeq 
A similar argument at $v_3$ gives 
$$
0=3h_4(t)+2\sum_{n\geq 1}h_4(t-2nl_4)
+2\sum_{n\geq 1}h_4(t-2nl_2)+2\sum_{n\geq 1}h_4(t-2nl_3)
$$
\beq
-2\sum_{n\geq 1}(g_1+f_2)(t-(2n-1)l_2)- 2\sum_{n\geq 1}g_1(t-(2n-1)l_3),\forall t.\label{key3}
\eeq 
By definition, we have $g_j(t)=0$ for $t<0$, and 
by unit speed of propagation and \rq{gh}, we have 
$h_j(t)=0$ for $t<\min (l_2,l_3),$ so the sums in \rq{key2},\rq{key3} are all finite.

To construct our shape control, we will use $f_2$ alone to attain the desired shape on $e_4$, while $f_1,f_2$ will be used jointly to attain control on the other edges. 
Recall we have  assumed $l_3\leq l_2;$ for presentational simplicity.
Recall
 $T_*=\max( l_1+l_2,l_3+l_4)$.
In the construction that follows, $f_2$ will be supported in the interval
$(T_*-\max(l_3+l_4,l_2),T_*))$ while $f_1$ will be supported in $(T_*-l_2-l_1, T_*). $ 

Let $\phi \in \mH^1$. For the moment, we assume $q=0$, pointing out the modifications necessary in the general case at the section's end. In what follows, we set $f_j(t)=g_k(t)=0$ for $t<0$ and all $j,k$, and $h_j(t)=0$ for $t<l_3$. When convenient, we will denote $u$ as $u^{(f_1,f_2)}.$

{\bf Step 1.} 

For this step,  we construct $f_2^1(t)$ supported 
in 
 $(T_*-l_4-l_3,T_*-l_3)$ so that 
$$u_4^{(0,f_2^1)}(x,T_*)=\phi_4(x).$$
Setting 
$v^{(h_4,0)}(x,T_*)=\phi_4(x)$
in \rq{nd}, we get $h_4(T_*-x)=-\phi_4'(x)$
on $(T_*-l_4,T_*)$, and  we then extend $h_4$ trivially to $(0,T_*).$
 We now compute the control $f_2^1$, uniquely determined  on 
 $(T_*-l_4-l_3,T_*-l_3)$, that generates $h_4$, i.e. such that
 $$u^{(0,f_2^1)}(v_3,t)=h_4(t), \ t\in (T_*-l_4,T_*).$$
Let 
$\Delta = \min (l_1,l_2-l_3,l_2) $ if $l_3>l_2$, and otherwise $\Delta =\min (l_1,l_2).$
Define  $t_n^*=n\Delta$ for $n\geq 0.$
We will use \rq{key2}, \rq{key3} to  solve  for  $g_1,f_2^1$ on successive intervals 
$[t_j^*,t_{j+1}^*]$ in terms of the known function $h_4$.
As $j$ increases, the number of non-zero terms in \rq{key2} and/or \rq{key3} increases. In what follows, we will denote by $\alpha (t)$ various functions, which change from line to line, that are known for $t<t_j^*$. Also, $\beta (t)$ are functions, that varying from line to line, that are determined by $h_4$, and hence are known for all $t$ in the argument below.

We  rewrite \rq{key2}, resp. \rq{key3} as: 
\beq
\al (t)=3g_1(t)+f_2^1(t)+2\sum_{n\geq 1}\mbox{ {$(-1)^n$}}g_1(t-2nl_1)
+2\sum_{n \geq 1}\big ( g_1(t-2nl_2)+f_2^1(t-2nl_2)\big )+2\sum_{n \geq 1}g_1(t-2nl_3),\label{key2'}
\eeq 
\beq
\be (t)=\sum_{n \geq 1}(g_1(t-(2n-1)l_2+l_3)+f_2(t-(2n-1)l_2+l_3)+g_1(t)+ \sum_{n \geq 1}g_1(t-2nl_3).\label{key3'}
\eeq
For $t<t_1^*$ and $n\geq 1$, we have $t-(2n-1)l_j<0$ for $j=1, 2,3$ by the definition of $t_j^*$, and hence
\beq 
\al(t)=3g_1(t)+f_2^1(t) ,\label{key2b}
\eeq
and similarly \rq{key3'} becomes
\beq
\be (t)=\big ( g_1(t-l_2+l_3)+f_2^1(t-l_2+l_3)\big )+ g_1(t).\label{key4}
\eeq
We have two cases.

Case 1:  $l_2>l_3.$
 Then for $t<t_1^*$, we have $t-l_2+l_3<0$, so
$$\beta (t)=g_1(t), \ t\in [0,t_{1}^*].$$
Combining this with \rq{key2b}, we can clearly solve for $f_2^1,g_1$ for $t<t_{1}^*.$

Case 2:  $l_2=l_3$.
Then \rq{key4} becomes
$$\alpha (t)=2g_1(t)+f_2^1(t), \ t\in [0,t_{1}^*], $$
and thus we solve for $f_2^1,g_1$ for $t<t_{1}^*.$


Now assume $f_2^1(t),g_1(t)$ have been determined for $t<t_k^*$. In this case, $t<t_{k+1}^*$ implies that for $n>0$ and $j=1,2,3$, we have 
$(t-2nl_j)<t_k^*$, and so the functions
$g_1(t-2nl_j),f_2^1(t-2nl_j)$ have been previously calculated. Thus \rq{key2'} becomes:
\beq
\al (t)=3g_1(t)+f_2^1(t),\ t\in [t_k^*,t_{k+1}^*] .\label{k22}
\eeq 
 
 We now rewrite \rq{key3'} by identifying known terms. First, if $t< t_{k+1}^*$ and $j=2,3$, then for $n>1$ we have $t-(2n-1)l_j+l_3<t_k^*,$ so  \rq{key4}  becomes
\beq
\al (t)=g_1(t-l_2+l_3)+f_2^1(t-l_2+l_3)+ 2g_1(t), \ t<t_{k+1}^*.\label{k32}
\eeq
By considering the same cases as above, we solve for $f_2^1(t),g_1(t)$ for $t<t_{k+1}^*$.






We iterate this argument until $t_{k+1}^*\geq T_*-l_3$; this determines  $f_2^1(t)$ for $t\in (T_*-l_4-l_3,T_*-l_3)$. We then extend $f_2^1$ trivially to the rest of $[0,T_*].$

{\bf Step 2.} 
Let $\psi=\phi -u^{(0,f_2^1)}(\cdot,T_*), $ so $\psi_4=0.$

In this step, we determine 
 $f_1^1$ supported on $(T_*-l_1-l_2,T_*-l_1)$, together with $f_2^2(t)$, supported on $(T_*-l_2,T_*),$
so that 
$$u^{(f_1^1,f_2^2)}_j(x,T_*)=\psi_j,\ j=2,3.$$
 Setting $f_2=f_2^1+f_2^2$, we will then have 
$$u^{(f_1^1,f_2)}_j(x,T_*)=\phi_j,\ j=2,3,4.$$

Since $q=0$, the integral kernels in  \rq{nd} vanish, so for $j=2,3$ we can easily find 
$g_j\in L^2(T_*-l_j,T_*)$ such that 
\beq
u^{g_j}(x,T_*)=v_{NN}^{(g_j,0)}(x,T_*)=\psi_j(x).\label{int}
\eeq 
We then set $g_3(t)=0$ for $t\leq  T_*-l_3,$ and $g_2(t)=0$ for $t\leq  T_*-l_2.$
Then $f_2^2$ is found immediately for all $t$  by applying \rq{g1g2},\rq{g1g3}, and clearly  $f_2^2$ has its support in $[T_*-l_2,T_*]$. The support of $f_2=f_2^1+f_2^2$ will  be 
$$[T_*-l_2,T_*]\cup [T_*-l_3-l_4,T_*-l_3]=
[T_*-\max (l_2,l_3+l_4),T_*].$$

We now solve for $(f_1^1)'$.
First note that by unit speed of wave propagation, for $t<T_*$ the waves generated by $g_2,g_3$ will not reach beyond $v_3, $ so $h_j(t)=0$. Thus, \rq{key2} simplifies to
$$
0=3g_1(t)+f_2^1(t)+2\sum_{n\geq 1} (-1)^ng_1(t-2nl_1)
+2\sum_{n\geq 1} (g_1(t-2nl_2)+f_2^1(t-2nl_2))+2\sum_{n\geq 1} g_1(t-2nl_3)
$$
\beq
+2\sum_{n\geq 1} (-1)^n(f_1^1)'(t-(2n-1)l_1), \ t\in (T_*-l_2,T_*).\label{key4'}
\eeq 
Since the functions $g_1,f_2^2$ are known, we can solve for $(f_1^1)'(t)$ iteratively for $t\in [T_*-l_2-l_1,T_*-l_1].$ We then set $f_1^1(t)=0$
for $t\leq T_*-l_2-l_1$.

{\bf Step 3.}
We have defined $f_2(t)$ for all $t$.
It remains to find $f_1(t)$. Choose some extension of  $f_1^1$ to $H^1_0(0,T_*).$
Let $\tilde{\psi}=\psi -u^{(f_1^1,f_2)}(*,T_*), $ so $\tilde{\psi}_j=0$ for $j=2,3,4.$
By continuity, we have $\tilde{\psi}(v_2)=0$, and also $\tilde{\psi}(v_1)=0$.
Then by  \rq{nd}, we solve $f_1^2\in H^1_0(T_*-l_1,T_*)$, with support in $(T_*-l_1,T_*)$,  so that $v_{ND}^{(f_1^2,0)}(x,T_*)=\tilde{\psi}_1(x),$ and hence
$u=u^{f_1^2,0}$ will solve $u_1(x,T_*)=\tilde{\psi}_1.$ We then  set $f_1=f_1^1+f_1^2$, and the shape control problem is solved.

We now indicate the changes necessary for the case
$q\neq 0$. In this case, the only two changes. First, solving \rq{int} now requires a standard integral equation argument. Second, 
 \rq{key4'} changes to 
$$
0=3g_1(t)+f_2(t)+2\sum_{n\geq 1} (-1)^ng_1(t-2nl_1)
+2\sum_{n\geq 1} (g_1(t-2nl_2)+f_2(t-2nl_2))+2\sum_{n\geq 1} g_1(t-2nl_3)
$$
\beq
+2\sum_{n\geq 1} (-1)^n(L_n(f_1^1)')(t-(2n-1)l_1), \ t<T_*,\label{key5}
\eeq 
with 
$$(L_n(f_1^1)')(t-(2n-1)l_1)=2(f_1^1)'(t-(2n-1)l_1)+2q((2n-1)l_1)f_1^1(t-(2n-1)l_1)$$
$$-\int_{s=0}^{t-(2n-1)l_1}\big ( k_+((2n-1)l_1,s+(2n-1)l_1)+k_-((2n-1)l_1,s+(2n-1)l_1)\big ) 
f_1^1(t-(2n-1)l_1-s)ds.$$
A standard argument using the theory of Volterra integral equations shows that $L_n$ is an isomorphism on $L^2(0,T_*-(2-1)l_1)$. Thus, using \rq{key5}, we solve for $L_n(f_1^1)'$, then $(f_1^1)'$, then $f_1^1$,  using at iterative argument. The rest of the proof of shape controllability is the same as in the case of $q=0$.$\Box$

{\bf Remark} 
Recalling that 	$T_*= \max \, \{\, l_1+l_2,\, l_3+l_4 \,\}$, we summarize the supports of the controls $f_1,f_2$.

1) if $l_1 + l_2 +\mu =  l_3 + l_4, \ \mu \geq 0,$  then $ T_* = l_3 + l_4, $ supp $\; f_1\subset [\mu,T_*]$ and  supp $\;f_2 \subset[ 0,T_* ]$; 

2)  if $l_1 + l_2=  l_3 + l_4 + \mu, \ \mu \geq 0,$  then $ T_* = l_1 + l_2,$ supp $\; f_1\subset [0,T_*]$ and  supp $\;f_2 \subset [ \tau,T_* ],$ 
where $\tau = \min\,\{\mu, l_1\}.$

\subsection{Proof of Theorem \ref{thm1}}
\label{proof}


Now we prove the exact controllability in time $2T,$ for $T\geq T_*,$
using the shape and velocity controllability in time $T$.
Let $\{ \om_n^2\}_1^{\infty}$ and  $\{ \f_n\}_1^{\infty}$   be the eigenvalues 
and  corresponding normalized eigenfunctions
of the spectral problem associated to system \rq{eqq1}-\rq{inc4}:
\beq \lab{sp}
- \f_n'' + q(x) \f_n = \om_n^2 \f_n, \ x\in e_j, j=1,..., 4, 
\eeq
\beq \lab{nc} 
(\f_n)_1 (v_1)=(\f_n)_4 (v_4)=0,
\eeq
\beq \lab{kn}
\pa (\f_n)_1(v_2)=\pa (\f_n)_2(v_2)=\pa (\f_n)_3(v_2), \ \sum_{j=1}^3  (\f_n)_j(v_2)=0,  
\eeq
\beq \lab{kn1}
\pa (\f_n)_2(v_3)=\pa (\f_n)_3(v_3)=\pa (\f_n)_4(v_3) , \ \sum_{j=2}^4  (\f_n)_j(v_3)=0. 
\eeq
Here and in what follows, we denote the restriction of $\f_n$ to edge $e_j$ by $(\f_n)_j.$

The following estimates can be extracted from  \cite{AN}, \cite{BK}, and \cite{Ku2}: 
$$|\pa (\f_{n})_1(v_1)|\prec  n, \ |(\f_n)_2(v_2)|\prec  1.
$$

We write
$u(x,t)=\sum_1^{\infty}a_n(t)\f_n(x). $ 
To find the coefficients $a_n$ we integrate by parts in $x$ the identity
$$0=\int_0^T \int_{\Omega} [\patt u-\paxx u+q(x)u]\,\f_n(x)\mu(t)\,dxdt$$
with arbitrary $\mu \in C^2_0[0,T].$ We obtain
the initial value problem
$$a_n''(t)+\omega_n^2a_n(t)=-f_1(t)
(\f_{n})_1'(v_1)
+f_{2}(t)(\f_n)_2(v_2), \ a_n(0)=a_n'(0)=0.
$$
By variation of parameters, we get
\beq \lab{an3}
a_n(T)=\int_0^T [-f_1(t)(\pa \f_{n})_1(v_1)
+f_2(t)(\f_n)_2(v_2)]\,\frac{\sin \om_n(T-t)}{\om_n}\,dt.
\eeq

In what follows, the controls solving the shape control problem will be denoted by 
$f_j=f_{j0}$, and the controls solving the velocity control problem will be denoted by $f_j=f_{j1},$ for $j=1,2.$
We apply integration by parts for the terms with $f_1,f_2$ in \rq{an3}, using $f_1(0)=f_1(T)=0$. Setting $a_n=a_n(T)$, the shape controllability result can be formulated as solvability of the moment problem
\beq \lab{3a}
{\om_n}a_n=\int_0^T f_{10}'(t)\ka_{n1}{\cos \om_n(T-t)}
+f_{20}(t)\ka_{n2}\,{\sin \om_n(T-t)}\,dt,\ \forall n, 
\eeq
for arbitrary sequence $\{a_n \om_n\} \in \ell^2$ by the functions $f_{10}',f_{20} \in L^2(0,T).$ Here $$\ka_{n1}=\frac{\pa (\f_{n})_1(v_1)}{\omn},\ \ka_{n2}={(\f_n)_{2}(v_2)}.$$
 
Differentiating \rq{an3} with respect to $T$, and  then integrating by parts while using $f_{1}(0)=f_{1}(T)=0$, we reformulate the 
the velocity controllability result  as solvability of the moment problem
\beq \lab{3bn}
b_n =\int_0^T [f_{21}(t)\,\ka_{n2}\,\cos \om_n(T-t)+ f'_{11}(t)\,\ka_{n1}\,\sin \om_n(T-t)]\,dt, \ \forall n,
\eeq
for arbitrary sequences $\{b_n\} \in \ell^2$ by the functions $f_{11}',f_{21} \in L^2(0,T).$ 

Now we extend the functions $f'_{11},f_{20}$ in the odd way with respect to $t=T$
from the interval $[0,T]$ to $[0,2T], $ and the functions $f'_{10},f_{21}$ --- in the even way. Define the functions 
$$ 
\ f'_{1}=\frac{f'_{10} + f'_{11}}{2}, \ 
f_{2}=\frac{f_{20}+f_{21}}{2},	$$
and observe that these functions solve the moment problems
\beq \lab{3a1}
a_n \om_n=\int_0^{2T} [f'_{1}(t)\,\ka_{n1}\,\cos \om_n(T-t)+f_{2}(t)\,\ka_{n2}\,\sin \om_n(T-t) ]\,dt
\eeq
\beq \lab{3b1}
b_n =\int_0^{2T} [f_{2}(t)\,\ka_{n2}\,\cos \om_n(T-t)+ f'_{1}(t)\,\ka_{n1}\,\sin \om_n(T-t)]\,dt
\eeq
for arbitrary sequences $\{a_n \om_n\}, \; \{b_n\} \in \ell^2$. Also, defining $f_{1}(t)=\int_0^t f'_{1}(s)\,ds,$  it's not hard to prove that $f_1\in H^1_0(0,2T)$.
Thus 
these moment problems can be rewritten in the form
\beq \lab{3a2}
a_n \om_n=\int_0^{2T} [-f_{2}(t)\,\ka_{n2}\ - f_{1}(t)\,\ka_{n1}\omega_n]\,\sin \om_n(T-t)\,dt
\eeq
\beq \lab{3b2}
b_n =\int_0^{2T} [-f_{2}(t)\,\ka_{n2}- f_{1}(t)\,\ka_{n1}\omega_n]
\,\cos \om_n(T-t)\,dt,
\eeq
 It is an algebra exercise to check that solvability of the moment problem \rq{3a2}, \rq{3b2} is
equivalent to solvability of the moment problem 
\beq \lab{3a3}
\al_n =\int_0^{2T} [-f_{2}(t)\,\ka_{n2}\ + f_{1}(t)\,\ka_{n1}\omega_n]\,\sin \om_n(2T-t)\,dt
\eeq
\beq \lab{3b3}
\be_n =\int_0^{2T} [-f_{2}(t)\,\ka_{n2}+ f_{1}(t)\,\ka_{n1}\omega_n]
\,\cos \om_n(2T-t)\,dt,
\eeq
for arbitrary sequences $\{\al_n\}, \{\be_n\} \in \ell^2.$
Having solved these moment problems is equivalent to proving exact controllability in time $2T$. The 
proof of Theorem \ref{thm1} is complete.$\Box$ 
		

\vskip3mm
\noindent {\bf  Acknowledgments}\\
The research of Sergei Avdonin was  supported  in part by the National Science Foundation,
	grant DMS 1909869.

\end{document}